\documentclass[12pt, leqno]{article}
\usepackage{amssymb,amsmath}
\usepackage[pdftex]{color}
\usepackage{latexsym}
\usepackage[pdftex]{graphicx}
\newtheorem{theorem}{Theorem}

\newtheorem{remark}{Remark}
\newtheorem{lemma}{Lemma}

\newcommand{\tx}[1]{\mbox{\;{#1}\;}} 
\newcommand{\N}{\mathbb{N}}
\newcommand{\R}{\mathbb{R}^n}

\newcommand{\Div}{\hbox{div}\:}
\newcommand{\p}{\partial}
\numberwithin{equation}{section}
\date{}
\usepackage[pdfproducer={LaTeX with hyperref package},pdfpagelayout=SinglePage,bookmarksnumbered,pdftex=true,breaklinks=true,bookmarks=true,linktocpage=true,pdfpagelabels=true,plainpages=false,pagebackref=false]{hyperref} 
\hypersetup{colorlinks=true,linkcolor=blue,citecolor=blue,urlcolor=blue,filecolor=blue}

\begin{document}
\title{Optimality conditions for the buckling of a clamped plate.}
\pagestyle{myheadings}
\maketitle
\centerline{\scshape Kathrin Knappmann}
\medskip
{\footnotesize
\centerline{Institut f\"ur Mathematik, RWTH Aachen  }
\centerline{Templergraben 55, D-52062 Aachen, Germany}}
\medskip
\centerline{\scshape Alfred Wagner}
\medskip
{\footnotesize
\centerline{Institut f\"ur Mathematik, RWTH Aachen  }
\centerline{Templergraben 55, D-52062 Aachen, Germany}}
%%%%%%%%%%%%%%%%%%%%%%%%%%%%%%%%%%%%%%%%%%%%%
\bigskip
\noindent
\abstract{We prove the following uniqueness result for the buckling plate. Assume there exists a smooth domain which minimizes the first buckling eigenvalue for a plate among all smooth domains of given volume. Then the domain must be a ball. The proof uses the second variation for the buckling eigenvalue and an inequality by L. E. Payne to establish this result.}
\bigskip

\bigskip

\centerline {MSC2010:  49K20, 49R05, 15A42, 35J20, 35N25.}
%%%%%%%%%%%%%%%%%%%%%%%%%%%%%%%%%%%%%%%%%%%%%
%%%%%%%%%%%%%%%%%%%%%%%%%%%%%%%%%%%%%%%%%%%%%
\section{Introduction}
We consider the following variational problem. Let $\Omega \subset \R$ be a bounded domain and let 
\begin{eqnarray*}
{\cal{R}}(u,\Omega):=\frac{\int\limits_{\Omega}\vert\Delta u\vert^2\:dx}{\int\limits_{\Omega}\vert\nabla u\vert^2\:dx}
\end{eqnarray*}
for $u\in H^{2,2}_0(\Omega)$. We set ${\cal{R}}(u,\Omega)=\infty$ if the denominator vanishes. We define
\begin{eqnarray}\label{eigen}
\Lambda(\Omega):=\inf\left\{ {\cal{R}}(u,\Omega)\: :\:u\in \: H^{2,2}_0(\Omega)\right\}.
\end{eqnarray}
The infimum is attained by the first eigenfunction $u$, which solves the Euler Lagrange equation
\begin{eqnarray}
\label{eq1} \Delta^2u+\Lambda(\Omega)\Delta u&=&0\quad\hbox{in}\:\Omega\\
\label{eq2} u=\partial_\nu u&=&0\quad\hbox{in}\:\partial\Omega.
\end{eqnarray}
If we normalize $u$ by $\| \nabla u\|_{L^{2}(\Omega)} =1$, the first eigenfunction is uniquely determined. Otherwise any multiple of $u$ is an eigenfunction as well. The sign of the first eigenfunction may change depending on $\Omega$.
\newline
\newline
The quantity $\Lambda(\Omega)$ is called buckling eigenvalue of $\Omega$.
It is well known that there is a discrete spectrum of positive eigenvalues of finite multiplicity and their only accumulation point is $\infty$. The corresponding eigenfunctions form an orthonormal basis of $H^{2,2}_0(\Omega)$. 
\newline
\newline
In the sequel, we will assume that $u$ is normalized. If we multiply \eqref{eq1} with $x\cdot\nabla u$  and integrate by parts, we obtain
\begin{eqnarray}\label{Rel}
\Lambda(\Omega) = \frac{1}{2}\int\limits_{\partial\Omega}\vert\Delta u\vert^2\, x\cdot\nu\:dS.
\end{eqnarray}
In 1951, G. Polya and G. Szeg\"o formulated the following conjecture (see \cite{PoSz51}). 
\newline\newline
\centerline{{\it{Among all domains $\Omega$ of given volume, the ball minimizes 
$\Lambda(\Omega).$}}} 
\newline\newline
This conjecture is still open. However, partial results are known. In \cite{Sz50} Szeg\"o 
proved the conjecture for all smooth domains under the additional assumption that 
$u > 0$ in $\Omega$. 
M.S. Ashbaugh and D. Bucur proved that among simply connected domains of 
prescribed volume there exists an optimal domain \cite{AB03}. 
In \cite{Wi95} H. Weinberg and B. Willms proved the following uniqueness result for $n=2$. 
If an optimal plane domain $\Omega$ exists and if $\partial\Omega$ is smooth 
(at least $C^{2,\alpha}$), then $\Omega$ is a disc. 
\newline\newline
There also exist bounds for $\Lambda(\Omega)$. 
We only mention Payne's inequality (see \cite{P55}) which states that 
\begin{eqnarray*}
\Lambda(\Omega) \geq \lambda_2(\Omega),
\end{eqnarray*}
where $\lambda_2$ denotes the second Dirichlet eigenvalue for the Laplacian. 
Equality holds if and only if $\Omega$ is a ball. 
\newline\newline
In this paper, we assume that there exists an optimal domain $\Omega$, which is smooth and simply connected. We will prove that $\Omega$ must be a ball. Thus we generalize the result of H. Weinberg and B. Willms in \cite{Wi95} to higher dimensions.
\newline\newline
To consider the second domain variation for $\Lambda(\Omega)$ is motivated by the work of E. Mohr in \cite{Mo75}. He was interested in the clamped plate eigenvalue, where 
\[
  {\cal{R}}(u,\Omega) = \frac{\int\limits_\Omega \vert\Delta u\vert^2\:dx}{\int\limits_\Omega u^2\:dx}
\]
and $\Omega$ is a smoothly bounded domain in $\mathbb{R}^2$. For the corresponding eigenvalue he computed the second domain variation. The explicit computation of the kernel of the second domain variation then implies that the disc is a unique minimizer among smooth domains of equal volume. 
\newline\newline
Our strategy will be as follows. In Chapter 2 we introduce a smooth family $(\Omega_t)_t$ of perturbations of $\Omega$ of equal volume. We denote by $\Lambda(t) := \Lambda(\Omega_t)$ the corresponding first buckling eigenvalue of $\Omega_t$. As a consequence of the optimality of $\Omega$, the eigenfunction $u$ statisfies the overdetermined boundary value problem
\begin{align*}
  \Delta^2 u + \Lambda(\Omega) \Delta u &=0 \; \mbox{in } \Omega \\
 u = \partial_\nu u &=0 \; \mbox{in } \partial\Omega \\
\Delta u &= c_0 \; \mbox{in } \partial\Omega, \mbox{ where } c_0 = \frac{2\Lambda(\Omega)}{\vert\Omega\vert} \text{ by } \eqref{Rel}.
\end{align*} 
This follows from the fact that the first domain variation of $\Lambda(\Omega)$ - computed in Chapter 3 - for an optimal domain necessarily vanishes.
\newline
In Chapter 4 we compute the second domain variation of $\Lambda(\Omega)$.  It turns out that 
\begin{equation}\label{ddotLambda}
  \ddot{\Lambda}(0) = \frac{d^2}{dt^2}\Lambda(t) \bigg{|}_{t=0} = 2\int\limits_\Omega\vert\Delta u'\vert^2 -2\Lambda(\Omega)\int\limits_\Omega\vert\nabla u'\vert^2\:dx,
\end{equation}
where $u'$ is the so called shape derivative of $u$. It solves 
\begin{eqnarray}
\label{shape1} \Delta^2 u' + \Lambda(\Omega) \Delta u' &=&0 \; \mbox{in } \Omega \\
\label{shape2} u'  &=&0 \; \mbox{in } \partial\Omega \\
\label{shape3}\partial_\nu u' &=& -c_0 v.\nu \; \mbox{in } \partial\Omega
\end{eqnarray} 
and 
\begin{eqnarray}
\label{shape4}\int\limits_\Omega \nabla u.\nabla u' \: dx =0. 
\end{eqnarray}
The vector field $v$ is the first order approximation of $\Omega_t$ in the sense that for $y \in \Omega_t$ there exists an $x\in\Omega$ such that 
\begin{eqnarray*}
y = x + tv(x) + o(t). 
\end{eqnarray*}
Thus, $\ddot{\Lambda}(0)$ is equal to a quadratic functional in the shape derivative $u'$ which we denote by ${\cal{E}}(u')$ and ${\cal{E}}(u')$ is given by the right hand side of \eqref{ddotLambda}. Since we assume the optimality of $\Omega$, we have ${\cal{E}}(u')\geq 0$. It turns out that the kernel of ${\cal{E}}(u')$ contains the directional derivatives $\partial_1 u,\ldots, \partial_n u$ of $u$. Each directional derivative is a shape derivative, which corresponds to a domain perturbation given by translations. 
\newline
\newline
The key idea is to enlarge the class of shape derivatives on which ${\cal{E}}$ is defined. This new class will be denoted by ${\cal{Z}}$ and contains the shape derivatives as a true subset. Nevertheless we can show that ${\cal{E}}$ is still bounded from below and even nonnegative on ${\cal{Z}}$. Moreover $\min_{{\cal{Z}}} {\cal{E}}=0$ since the directional derivatives of $u$ are in ${\cal{Z}}$. This is done in Chapter 5. In Chapter 6 we construct a function $\psi\in{\cal{Z}}$ for which we will show
\begin{eqnarray*}
0\leq {\cal{E}}(\psi)\leq \left(\lambda_2(\Omega)-\Lambda(\Omega)\right)\lambda_2(\Omega).
\end{eqnarray*}
By Payne's inequality we have equality and this proves that the optimal domain is a ball.
\newline\newline
Some of these results were obtained in the Diplom thesis of the first author \cite{Kn08}.
%%%%%%%%%%%%%%%%%%%%%%%%%%%%%%%%%%%%%%%%%%%%%
%%%%%%%%%%%%%%%%%%%%%%%%%%%%%%%%%%%%%%%%%%%%%
\section{Domain variations}\label{DomVar}
%%%%%%%%%%%%%%%%%%%%%%%%%%%%%%%%%%%%%%%%%%%%%
%%%%%%%%%%%%%%%%%%%%%%%%%%%%%%%%%%%%%%%%%%%%%
Let $\Omega$ be a bounded smooth (at least $C^{2,\alpha}$) and simply connected domain in $\R$. We denote by $\nu$ the unit normal vector field on $\partial\Omega$. Let $\delta$ be the distance function to the boundary, i.e. for $x\in\overline{\Omega}$ we have 
\begin{eqnarray*}
\delta(x):=\inf\{\vert x-z\vert\: :\:z\in\,\partial\Omega\}.
\end{eqnarray*} 
Then, for smooth $\partial\Omega$, $\nu:=\nabla\delta$ defines a smooth extension of $\nu$ into a sufficiently small tubular neighbourhood of $\partial\Omega$. With this the following identities hold.
\begin{eqnarray}\label{normal}
\nu\cdot\nu=1,\qquad \nu\cdot D\nu=0\qquad\hbox{and}\qquad D\nu\cdot\nu=0
\end{eqnarray}
on $\partial\Omega$. See e.g. Proposition 5.4.14 in \cite{HePi05} for a proof.
\newline
\newline
Moreover, the mean curvature of $\partial\Omega$ is bounded since $\Omega$ is smooth, i.e. for each $x \in \partial\Omega$ there holds
\begin{equation}\label{H_dOmega}
  \vert H_{\partial\Omega}(x) \vert \leq \max_{\partial\Omega} \vert H_{\partial\Omega}\vert < \infty.
\end{equation}
We will frequently use integration by parts on $\partial\Omega$. Let $f\in C^1(\partial\Omega)$ and $v\in C^{0,1}(\partial\Omega,\R)$. The next formula is often called the {\sl Gauss theorem on surfaces}.
\begin{eqnarray}\label{tanpi}
\oint_{\p\Omega}f\:\Div_{\partial\Omega}v\:dS
=-\oint_{\p\Omega}v\cdot\nabla^\tau f\:dS
+(n-1)\oint_{\p\Omega}f(v\cdot\nu)\:H_{\partial\Omega}\:dS,
\end{eqnarray}
where
\begin{align}\label{tgradient}
\nabla^{\tau} f =\nabla f -(\nabla f\cdot \nu)\nu
\end{align}
denotes the tangential gradient of $f$.
\newline
\newline
In this chapter, we describe the class of admissible variations for the domain functional 
$\Lambda(\Omega)$. 
For given $t_0>0$ and  $t\in(-t_0,t_0)$ let $(\Omega_t)_t$ be a family of perturbations of the 
domain $\Omega\subset \R$ of the form
\begin{eqnarray*}
\Omega_t=\Phi_{t}(\Omega)
\end{eqnarray*}
where
\begin{eqnarray*}
\Phi_{t}\::\:\overline{\Omega}\to\R
\end{eqnarray*}
is a diffeomorphism which is smooth in $t$ and $x$. Thus we may write
\begin{eqnarray*}
\Omega_t:=\{y=x+tv(x)+\frac{t^2}{2}w(x)+o(t^2):x \in \Omega,\: t\tx{small}\},
\end{eqnarray*}
where 
\begin{eqnarray*}
v=(v_1(x),v_2(x),\dots,v_n(x))=\partial_{t}\Phi_{t}(x)\vert_{t=0}
\end{eqnarray*}
and 
\begin{eqnarray*}
w=(w_1(x),w_2(x),\dots,w_n(x))=\partial^2_{t}\Phi_{t}(x)\vert_{t=0}
\end{eqnarray*} 
are smooth
vector fields and where $o(t^2)$ collects terms such that $\frac{o(t^2)}{t^2}\to 0$ as $t\to 0$. 
For small $t_0$ the sets $\Omega_t$ and $\Omega$ are diffeomorphic. We will frequently use the notation $y:=\Phi_{t}(x)$. Consider the functional 
\begin{eqnarray*}%\label{t-eigenvalue}
\Lambda(\Omega_t):=\inf\left\{ {\cal{R}}(u,\Omega_t)\: :\:u\in \: H^{2,2}_0(\Omega_t)\right\},
\end{eqnarray*} 
which only depends on $\Omega_t$. Let $u(t,y)\in H^{2,2}_0(\Omega_t)$ be the minimizer. For short we will write
\begin{eqnarray}\label{t-u}
\tilde{u}(t):=u(t,y).
\end{eqnarray}
Then $\tilde{u}(t)$ solves
\begin{eqnarray}
\label{t-eq1}\Delta^2\tilde{u}(t)+\Lambda(\Omega_t)\Delta\tilde{u}(t)
&=&0\quad\hbox{in}\:\Omega_t\\
\label{t-eq2}\tilde{u}(t)=\vert\nabla\tilde{u}(t)\vert&=&0\quad\hbox{in}\:\partial\Omega_t.
\end{eqnarray}
for each $t\in(-t_0,t_0)$.
With this notation we define
\begin{eqnarray*}%\label{t-lambda}
\Lambda(t):= {\cal{R}}(\tilde{u}(t),\Omega_t).
\end{eqnarray*}
Since we assume smoothness of $\Omega$ and $\Phi_t$ the eigenfunction $\tilde{u}$ is also smooth in $t$ and $x$. This has several consequences which we list as remarks.
%%%%%%%%%%%%%%%%%%%%%%%%%%%%%%%%%%%%%%%%
\begin{remark}
Since $\partial\Omega_t$ is smooth and since $\tilde{u}(t)=0$ on $\partial\Omega_t$ then necessarily
\begin{eqnarray}
\label{lapbc1}\Delta \tilde{u}=\partial_{\nu}^2\tilde{u}+(n-1)\partial_{\nu}\tilde{u} \:H_{\partial\Omega_t}\quad\hbox{in}\:\partial\Omega_t,
\end{eqnarray}
where $H_{\partial\Omega_t}$ denotes the mean curvature of $\partial\Omega_t$. Clearly,
if $\tilde{u}=\vert\nabla \tilde{u}\vert =0$ on $\partial\Omega_t$, then necessarily
\begin{eqnarray}
\label{lapbc2}\Delta \tilde{u}=\partial_{\nu}^2\tilde{u}\quad\hbox{in}\:\partial\Omega_t.
\end{eqnarray}
\end{remark}
%%%%%%%%%%%%%%%%%%%%%%%%%%%%%%%%%%%%%%%%
\begin{remark}
Since \eqref{t-eq2} holds for all $t\in(-t_0,t_0)$, we also have
\begin{eqnarray}
\label{tdot-eq2}\dot{\tilde{u}}(t)=\vert\nabla\dot{\tilde{u}}(t)\vert=0\quad\hbox{in}\:\partial\Omega_t
\end{eqnarray}
for all $t\in(-t_0,t_0)$.
\end{remark}
%%%%%%%%%%%%%%%%%%%%%%%%%%%%%%%%%%%%%%%%
\begin{remark}
Straightforward computations yield
\begin{eqnarray*}
\dot{\tilde{u}}(t)=\frac{d}{dt}u(t,y)=\partial_{t}u(t,\Phi_{t}(\Phi^{-1}_{t}(y))+\partial_{t}\Phi_{t}(\Phi^{-1}_{t}(y))\cdot\nabla u(t,y)
\end{eqnarray*}
for all $t\in(-t_0,t_0)$. Let $y\in\partial\Omega_t$. Then \eqref{tdot-eq2}  and \eqref{t-eq2} imply
\begin{eqnarray}
\label{tdot-eq3} 0=\dot{\tilde{u}}(t)=\partial_{t}u(t,y)\quad\hbox{for $y$ in}\:\partial\Omega_t
\end{eqnarray}
for all $t\in(-t_0,t_0)$.
\end{remark}
%%%%%%%%%%%%%%%%%%%%%%%%%%%%%%%%%%%%%%%%
In particular for $t=0$ we compute $\tilde{u}(0)=u(x)$ and
\begin{eqnarray*}
\dot{\tilde{u}}(0)&=& \partial_{t}u(0,x)+v(x)\cdot Du(0,x)\\
\ddot{\tilde{u}}(0)&=& \partial^2_{t}u(0,x)+2v(x)\cdot D\partial_{t}u(0,x)+w(x)\cdot Du(0,x)
+v(x)\cdot D\left(v(x)\cdot Du(0,x)\right).
\end{eqnarray*}
We will use the notation
\begin{eqnarray*}
u'(x):=\partial_{t}u(0,x)\qquad\hbox{and}\qquad u''(x):=\partial^2_{t}u(0,x).
\end{eqnarray*}
Hence,
\begin{eqnarray}
\label{tildeu1} \dot{\tilde{u}}(0)&=& u'(x)+v(x)\cdot Du(x)\\
\label{tildeu2} \ddot{\tilde{u}}(0)&=& u''(x)+2v(x)\cdot Du'(x)+w(x)\cdot Du(x)
+v(x)\cdot D\left(v(x)\cdot Du(x)\right).
\end{eqnarray}
Note that all these quantities are defined for $x\in\overline{\Omega}$. For $x\in\partial\Omega$ we thus get
\begin{eqnarray*}
0 = \dot{\tilde{u}}(0) = u'(x)\qquad \hbox{and} \qquad  0=\nabla\dot{\tilde{u}}(0) =\nabla u'(x)+v(x)\cdot D^2u(x),
\end{eqnarray*}
where $(v(x)\cdot D^2u(x))_{j}=\sum_{i=1}^{n}v_{i}(x)\partial_{i}\partial_{j}u(x)$ for $j=1,\hdots,n$. Thus, we get the following boundary conditions for $u'$.
\begin{eqnarray}\label{uprimebc}
u'(x)=0\qquad\hbox{and}\qquad \partial_{\nu}u'(x)=-v(x)\cdot D^2u(x)\cdot\nu(x)\quad \hbox{for}\quad x\in\partial\Omega.
\end{eqnarray}
Here we used the notation $v(x)\cdot D^2u(x)\cdot\nu(x)=\sum_{i,j=1}^{n}v_{i}(x)
\partial_{i}\partial_{j}u(x)\nu_{j}(x)$.
\newline
\newline
%%%%%%%%%%%%%%%%%%%%%%%%%%%%%%%%%%%%%%%%
Let $\nu_{t}(y)$ be the unit normal vector in $y\in\partial\Omega_t$. We also write this as
\begin{eqnarray}
\label{normal1}\nu_{t}(y)=\nu(t,\Phi_{t}(x))\qquad \forall\:t\in(-t_0,t_0)\qquad x\in\partial\Omega.
\end{eqnarray}
Then we have
\begin{eqnarray}\label{normal2}
\nu'=-\nabla^{\tau}(v\cdot\nu),\qquad\nu\cdot\nu'=0.
\end{eqnarray}
This follows from direct calculations (see e.g. (5.64) in \cite{HePi05}).
%%%%%%%%%%%%%%%%%%%%%%%%%%%%%%%%%%%%%%%%
\begin{lemma}
With the notation from above the following equality holds.
\begin{eqnarray}
\label{tdot-eq4} \nu_{t}\cdot\nabla(\partial_{t}u(t,y))=-\Delta u(t,y)\:\nu_{t}\cdot \partial_{t}\Phi_{t}(\Phi^{-1}_{t}(y))\quad\hbox{for $y$ in}\:\partial\Omega_t
\end{eqnarray}
for all $t\in(-t_0,t_0)$. Alternatively, we write this for all $t\in(-t_0,t_0)$ and $x\in\partial\Omega$ as
\begin{eqnarray}
\label{tdot-eq5} \nu(t,\Phi_{t}(x))\cdot\nabla\{\partial_{t}u(t,\Phi_{t}(x))\}=-\Delta u(t,\Phi_{t}(x))\:\nu({t},\Phi_{t}(x))\cdot \partial_{t}\Phi_{t}(x).
\end{eqnarray}
\end{lemma}
%%%%%%%%%%%%%%%%%%%%%%%%%%%%%%%%%%%%%%%%
{\bf{Proof}} Since $\nabla u(t,\Phi_{t}(x))=0$ for all $\vert t\vert<t_0$ and all $x\in\partial\Omega$ we have
\begin{eqnarray*}
0=\frac{d}{dt}\:\nabla u(t,\Phi_{t}(x))=\nabla \partial_{t}u(t,\Phi_{t}(x)) + D^2u(t,\Phi_{t}(x))\cdot \partial_{t}\Phi_{t}(x).
\end{eqnarray*}
This implies
\begin{eqnarray*}
0=\nu_{t}\cdot\nabla(\partial_{t}u(t,y))+\nu_{t}\cdot D^2u(t,y)\cdot \partial_{t}\Phi_{t}(\Phi^{-1}_{t}(y))\quad\hbox{for $y$ in}\:\partial\Omega_t
\end{eqnarray*}
for all $t\in(-t_0,t_0)$. Here we used the notation
\begin{eqnarray*}
\nu_{t}\cdot D^2u(t,y)\cdot \partial_{t}\Phi_{t}(\Phi^{-1}_{t}(y))
=
\sum_{ij=1}^{n}\nu_{t,i}\cdot \partial_{i}\partial_{j}u(t,y)\cdot \partial_{t}\Phi_{t}(\Phi^{-1}_{t}(y))_{j}.
\end{eqnarray*}
Since $\nabla\tilde{u}(t)=0$ in $\partial\Omega_t$, we get
\begin{eqnarray*}
\nu_{t}\cdot D^2u(t,y)\cdot \partial_{t}\Phi_{t}(\Phi^{-1}_{t}(y))
=
\nu_{t}\cdot D^2u(t,y)\cdot\nu_{t}\:\nu_{t}\cdot \partial_{t}\Phi_{t}(\Phi^{-1}_{t}(y)).
\end{eqnarray*}
Thus,
\begin{eqnarray*}
\nu_{t}\cdot\nabla(\partial_{t}u(t,y))=-\nu_{t}\cdot D^2u(t,y)\cdot\nu_{t}\:\nu_{t}\cdot \partial_{t}\Phi_{t}(\Phi^{-1}_{t}(y))\quad\hbox{for $y$ in}\:\partial\Omega_t.
\end{eqnarray*}
Formula \eqref{lapbc2} simplifies to
\begin{eqnarray*}
\nu_{t}\cdot\nabla(\partial_{t}u(t,y))=-\Delta u(t,y)\:\nu_{t}\cdot \partial_{t}\Phi_{t}(\Phi^{-1}_{t}(y))\quad\hbox{for $y$ in}\:\partial\Omega_t.
\end{eqnarray*}
This proves the lemma.
\hfill $\square$
\newline
\newline
%%%%%%%%%%%%%%%%%%%%%%%%%%%%%%%%%%%%%%%%
The first derivative of $\Lambda(t)$ with respect to the parameter $t$  is called 
{\sl the first domain variation} and the second derivative is called 
{\sl  the second domain variation}. 
\newline
\newline
Our domain variations will be chosen within the class of volume preserving pertubations up to
 order $2$. Hence, they are chosen such that
\begin{eqnarray}\label{volume}
{\cal{L}}^n(\Omega_t)={\cal{L}}^n(\Omega) +o(t^2)
\end{eqnarray}
holds. This puts constraints on the vector fields $v$ and $w$. They were discussed e.g. in 
\cite{BaWa14}, formula (2.13) and Lemma 1. 
%%%%%%%%%%%%%%%%%%%%%%%%%%%%%%%%%%%%%%%%%%%%%
\begin{lemma}\label{volpres1} Let $v,w\in C^{0,1}(\Omega,\R)$ be such that 
\eqref{volume} holds. Then
\begin{align}\label{volume1}
\int\limits_{\Omega} \Div v\:dx =0
\end{align}
and
\begin{eqnarray*}
\int\limits_{\Omega}\left((\Div v)^2-Dv : Dv+\Div\:w\right)\:dx=0,
\end{eqnarray*}
where $Dv:Dv=\sum_{i,j=1}^{n}\partial_{i}v_{j}\:\partial_{j}v_{i}$.
The second equality is equivalent to
\begin{align}\label{volume2}
\int\limits_{\partial\Omega}(v\cdot\nu) \Div v\:dS-\int\limits_{\partial\Omega}v\cdot Dv\cdot\nu\:dS
+
\int\limits_{\partial\Omega}(w\cdot \nu)\:dS=0.
\end{align}
\end{lemma}
%%%%%%%%%%%%%%%%%%%%%%%%%%%%%%%%%%%%%%%%
Note that rotations do not satisfy these conditions (see e.g. Remark 1 in \cite{BaWa14}).
%%%%%%%%%%%%%%%%%%%%%%%%%%%%%%%%%%%%%%%%
%%%%%%%%%%%%%%%%%%%%%%%%%%%%%%%%%%%%%%%%
\section{The first domain variation}
%%%%%%%%%%%%%%%%%%%%%%%%%%%%%%%%%%%%%%%%
%%%%%%%%%%%%%%%%%%%%%%%%%%%%%%%%%%%%%%%%
We will use the following formula for the computations of the first domain variation of 
$\Lambda$. It is well known as Reynolds transport theorem and is analyzed in detail 
in Chapter 5.2.3 in \cite{HePi05}. 
%%%%%%%%%%%%%%%%%%%%%%%%%%%%%%%%%%%%%%%%%%%
\begin{theorem}\label{reynolds1}
Let $t\in(-t_0,t_0)$ for some $t_0>0$. Let $\Phi_{t}\in C^{0,1}(\R)$ be differentiable in $t$ and let $t\to f(t)\in L^1(\R)$ be a function which is differentiable in $t$. Moreover, let $f(t)\in W^{1,1}(\R)$. Then $t\to I(t):=\int\limits_{\Omega_t}f(t)\:dy$ is differentiable in $t$. Moreover, we have the formula
\begin{eqnarray*}
\dot{I}(t)=\int\limits_{\Omega_t}\partial_{t}f(t)+\Div\left(f(t)\: \partial_{t}\Phi_{t}(\Phi^{-1}_{t}(y))\right)\:dy.
\end{eqnarray*}
If $\partial\Omega$ is sufficiently smooth (at least Lipschitz continuous), this is equivalent to
\begin{eqnarray*}
\dot{I}(t)=\int\limits_{\Omega_t}\partial_{t}{f}(t)\:dy+\int\limits_{\partial\Omega_t}f(t)\: \partial_{t}\Phi_{t}(\Phi^{-1}_{t}(y))\cdot\nu(y)\:dS(y).
\end{eqnarray*}
In particular, for $t=0$ we get
\begin{eqnarray*}
\dot{I}(0)=\int\limits_{\Omega}\partial_{t} f (t)\vert_{t=0}+\Div\left(f(0)\: v(x)\right)\:dx.
\end{eqnarray*}
Again, if $\partial\Omega$ is sufficiently smooth, this is equivalent to
\begin{eqnarray*}
\dot{I}(0)=\int\limits_{\Omega}\partial_{t} f (t)\vert_{t=0}\:dx+\int\limits_{\partial\Omega}f(0)\: v(x)\cdot\nu(x)\:dS(x).
\end{eqnarray*}
\end{theorem}
%%%%%%%%%%%%%%%%%%%%%%%%%%%%%%%%%%%%%%%%%%%
We apply this formula to $\Lambda(t)=\frac{D(t)}{N(t)}$ where
\begin{eqnarray*}
D(t):=\int\limits_{\Omega_t}\vert\Delta\tilde{u}(t)\vert^2\:dy\qquad\hbox{and}\qquad
N(t):=\int\limits_{\Omega_t}\vert\nabla\tilde{u}(t)\vert^2\:dy
\end{eqnarray*}
and we assume the normalization
\begin{eqnarray}\label{norm}
N(t)=\int\limits_{\Omega_t}\vert\nabla\tilde{u}(t)\vert^2\:dy=1\quad\forall\:t\in(-t_0,t_0).
\end{eqnarray} 
We then obtain
\begin{eqnarray*}
\label{ladot}\dot{\Lambda}(t)&=&2\int\limits_{\Omega_t}\Delta\tilde{u}(t)\:\Delta\partial_{t}{\tilde{u}}(t)\:dy
-
2\:\Lambda(t)\:\int\limits_{\Omega_t}\nabla\tilde{u}(t)\cdot\nabla\partial_{t}{\tilde{u}}(t)\:dy\\
&&+
\nonumber\int\limits_{\partial\Omega_t}\vert\Delta\tilde{u}(t)\vert^2\:\partial_{t}\Phi_{t}(\Phi_{t}^{-1}(y))\cdot\nu_{t}(y)\:dS(y),
\end{eqnarray*} 
where $\nu_{t}(y)$ denotes the unit normal vector in $y\in\partial\Omega_t$. 
We integrate by parts and use \eqref{tdot-eq2}. Then
\begin{eqnarray*}
\dot{\Lambda}(t)&=&
2\int\limits_{\Omega_t}\left\{\Delta^2\tilde{u}(t)+\Lambda(t)\:\Delta\tilde{u}(t)\right\}\partial_{t}{\tilde{u}}(t)\:dy
+
2\int\limits_{\partial\Omega_t}\Delta\tilde{u}(t)\:\partial_{\nu_t}\partial_{t}{\tilde{u}}(t)\:dS(y)\\
&&-
2\int\limits_{\partial\Omega_t}\partial_{\nu_t}\Delta\tilde{u}(t)\partial_{t}{\tilde{u}}(t)\:dS(y)
+
\int\limits_{\partial\Omega_t}\vert\Delta\tilde{u}(t)\vert^2\:\partial_{t}\Phi_{t}(\Phi_{t}^{-1}(y))\cdot\nu_{t}(y)\:dS(y).
\end{eqnarray*} 
The first integral vanishes since $\tilde{u}(t)$ solves \eqref{t-eq1}. The third integral vanishes since \eqref{tdot-eq3} holds. Finally we use \eqref{tdot-eq4}. This proves the following lemma.
%%%%%%%%%%%%%%%%%%%%%%%%%%%%%%%%%%%%%%%%%
\begin{lemma}\label{ladot1} 
Let $\tilde{u}(t)$ be an eigenfunction (i.e. a solution of \eqref{t-eq1} - \eqref{t-eq2}) and assume \eqref{norm} holds. Let
\begin{eqnarray*}
\Lambda(t)=\int\limits_{\Omega_t}\vert\Delta\tilde{u}(t)\vert^2\:dy.
\end{eqnarray*}
Then
\begin{eqnarray}\label{ladotform}
\dot{\Lambda}(t)&=&
-
\int\limits_{\partial\Omega_t}\vert\Delta\tilde{u}(t)\vert^2\:\partial_{t}\Phi_{t}(\Phi_{t}^{-1}(y))\cdot\nu_{t}(y)\:dS(y).
\end{eqnarray}
\end{lemma}
%%%%%%%%%%%%%%%%%%%%%%%%%%%%%%%%%%%%%%%%%
\begin{remark}\label{mon}
Note that if $\partial_{t}\Phi_{t}(\Phi_{t}^{-1}(y))\cdot\nu_{t}(y)>0$, this implies 
${\cal{L}}^n(\Omega_t)>{\cal{L}}^n(\Omega)$ for small $t$. Thus, $\dot{\Lambda}(t)$ is negative in this case. From this we conclude that the first buckling eigenvalue is decreasing under set inclusion.
\end{remark}
%%%%%%%%%%%%%%%%%%%%%%%%%%%%%%%%%%%%%%%%%
From Lemma \ref{ladot1} we get in particular
\begin{eqnarray*}
\dot{\Lambda}(0)=-\int\limits_{\partial\Omega}\vert\Delta u\vert^2\:v(x)\cdot\nu(x)\:dS(x).
\end{eqnarray*}
From Lemma \ref{volpres1} and \eqref{volume1} we deduce $\vert\Delta u\vert=const.$ if $\Omega$ is a critical point of $\Lambda(t)$. 
Due to formula \eqref{Rel}, this constant is equal to 
\begin{eqnarray}\label{constant}
c_0 := \frac{2\Lambda(0)}{\vert\Omega\vert}.
\end{eqnarray}
We denote this result as a theorem.
%%%%%%%%%%%%%%%%%%%%%%%%%%%%%%%%%%%%%%%%%
\begin{theorem}\label{overdet1}
Let $\Omega_t$ be a family of volume preserving perturbations  of $\Omega$ as described in Chapter 2. Then $\Omega$ is a critical point of the energy $\Lambda(t)$, i.e. $\dot{\Lambda}(0)=0$,
if and only if
\begin{eqnarray}\label{overrt}
\Delta u=c_0\qquad\hbox{on}\quad\partial\Omega.
\end{eqnarray}
\end{theorem}
%%%%%%%%%%%%%%%%%%%%%%%%%%%%%%%%%%%%%%%%%
In particular, $u$ is a solution of the overdetermined boundary value problem
\begin{eqnarray}
\label{ceq1} \Delta^2u+\Lambda(\Omega)\Delta u&=&0\quad\hbox{in}\:\Omega\\
\label{ceq2} u=\partial_{\nu}\nabla u&=&0\quad\hbox{in}\:\partial\Omega.\\
\label{ceq3} \Delta u&=&c_0>0\quad\hbox{in}\:\partial\Omega.
\end{eqnarray}
Note that if we set $U:=\Delta u+\Lambda(\Omega) u$ \eqref{ceq1} - \eqref{ceq3} implies
\begin{eqnarray*}
\Delta U=0\:\hbox{in}\: \Omega\: \hbox{ and } \: U=c_0\:\hbox{in}\: \partial\Omega.
\end{eqnarray*}
Hence,
\begin{eqnarray}\label{w-eq}
U= \Delta u+\Lambda(\Omega) u=c_0\qquad\hbox{in}\:\:\overline{\Omega}.
\end{eqnarray}
From \cite{Wi95} we know that for $n=2$ this implies that $\Omega$ is a ball.
In particular,
\begin{eqnarray}\label{w-eq1}
\partial_{\nu}\Delta u=0\qquad\hbox{in}\:\:\partial{\Omega}.
\end{eqnarray}
%%%%%%%%%%%%%%%%%%%%%%%%%%%%%%%%%%%%%%%%%
%%%%%%%%%%%%%%%%%%%%%%%%%%%%%%%%%%%%%%%%%
\section{The second domain variation}
%%%%%%%%%%%%%%%%%%%%%%%%%%%%%%%%%%%%%%%%%
%%%%%%%%%%%%%%%%%%%%%%%%%%%%%%%%%%%%%%%%%
Throughout this chapter we assume that $\Omega$ is an optimal domain, i.e. $\dot{\Lambda}(0)=0$ and $\ddot{\Lambda}(0)\geq 0$. This implies that $u$ solves \eqref{ceq1} -  \eqref{ceq3} and \eqref{w-eq}. As a consequence \eqref{uprimebc} reads as
\begin{eqnarray}
\label{upricc} u'(x)=0\qquad\hbox{and}\qquad\partial_{\nu}u'(x)=-c_0\:v(x)\cdot \nu(x)\quad \hbox{for}\quad x\in\partial\Omega.
\end{eqnarray}
Note that if we differentiate \eqref{t-eq1} - \eqref{t-eq2} in $t=0$ and use the fact that $\dot{\Lambda}(0)=0$, we obtain an equation for $u'$:
\begin{eqnarray}\label{uprimeq}
\Delta^2 u'(x)+\Lambda(\Omega)\Delta u'(x)=0\qquad\hbox{in}\:\:\Omega.
\end{eqnarray}
The boundary conditions for $u'$ are given by \eqref{upricc}.
Furthermore, the normalization \eqref{norm} implies 
\begin{eqnarray}\label{graduprime}
\int\limits_\Omega \nabla u\cdot\nabla u'\:dx =0.
\end{eqnarray}
We recall formula \eqref{ladotform}. Before we differentiate with respect to $t$ again we state the following consequence of Reynold's theorem (see e.g. Chapter 5.4.2 in \cite{HePi05}).
%%%%%%%%%%%%%%%%%%%%%%%%%%%%%%%%%%%%%%%%%
\begin{theorem}\label{reynolds2}
Let $\Omega$ be a bounded smooth domain of class $C^3$. Let $t\in(-t_0,t_0)$ and let $\Phi_{t}\in C^{0,1}(\R)$ be differentiable in $t$. Let $t\to g(t)\in L^1(\R)$ be a function which is differentiable in $t$. Moreover, let $g(t)\in W^{1,1}(\R)$. Then $t\to J(t):=\int\limits_{\partial\Omega_t}g(t)\:dS(y)$ is differentiable in $t$. Moreover, for $t=0$ we have the formula
\begin{eqnarray*}
\dot{J}(0)=\int\limits_{\partial\Omega}\partial_{t}g(0)+(v(x)\cdot\nu(x))\left\{\partial_{\nu}g(0)+(n-1)g(0)\:H_{\partial\Omega}(x)\right\}\:dS(x),
\end{eqnarray*}
where $H_{\partial\Omega}$ denotes the mean curvature of $\partial\Omega$ in $x$.
\end{theorem}
%%%%%%%%%%%%%%%%%%%%%%%%%%%%%%%%%%%%%%%%%
We apply this theorem to \eqref{ladotform}. It is convenient to apply \eqref{tdot-eq4} and to rewrite \eqref{ladotform} as
\begin{eqnarray*}
\dot{\Lambda}(t)&=&
\int\limits_{\partial\Omega_t}\Delta\tilde{u}(t)\:\nu_{t}\cdot\nabla(\partial_{t}u(t,y))\:dS(y).
\end{eqnarray*}
Let
\begin{eqnarray*}
g(t):=\Delta\tilde{u}(t)\:\nu_{t}\cdot\nabla(\partial_{t}u(t,y)).
\end{eqnarray*}
An application of Theorem \ref{reynolds2} yields
\begin{eqnarray}
\ddot{\Lambda}(0)&=&
\label{ddlamb}\int\limits_{\partial\Omega}\Delta u'\:\partial_{\nu}u'\:dS
+
\int\limits_{\partial\Omega}\Delta u\:\nu'\cdot\nabla u'\:dS
+
\int\limits_{\partial\Omega}\Delta u\:\partial_{\nu}u''\:dS\\
\nonumber&&+
\int\limits_{\partial\Omega}(v\cdot\nu)\:\partial_{\nu}(\Delta u\:\partial_{\nu}u')\:dS
+(n-1)
\int\limits_{\partial\Omega}(v\cdot\nu)\:\Delta u\:\partial_{\nu}u'\:H_{\partial\Omega}\:dS.
\end{eqnarray}
Note that 
\begin{eqnarray*}
\nu_{t}\cdot\nu_{t}=1\:\:\hbox{in}\:\:\partial\Omega_{t}\quad\Longrightarrow\quad\nu\cdot\nu'=0\:\:\hbox{in}\:\:\partial\Omega,
\end{eqnarray*}
where 
\begin{eqnarray*}
\nu'(x)=\partial_{t}\nu(t,\Phi_{t}(x))\vert_{t=0}\quad \hbox{for} \quad x\in\partial\Omega.
\end{eqnarray*}
Since \eqref{upricc} implies $\nabla u'=\partial_{\nu}u'\:\nu$, this implies
\begin{eqnarray*}
\int\limits_{\partial\Omega}\Delta u\:\nu'\cdot\nabla u'\:dS=0.
\end{eqnarray*}
For the fourth integral we apply \eqref{overrt} and \eqref{w-eq1}.
\begin{eqnarray*}
\int\limits_{\partial\Omega}(v\cdot\nu)\:\partial_{\nu}(\Delta u\:\partial_{\nu}u')\:dS
&=&
\int\limits_{\partial\Omega}(v\cdot\nu)\:\partial_{\nu}\Delta u\:\partial_{\nu}u'\:dS
+
\int\limits_{\partial\Omega}(v\cdot\nu)\:\Delta u\:\partial^2_{\nu}u'\:dS\\
\nonumber&&=
0+c_0\:\int\limits_{\partial\Omega}(v\cdot\nu)\:\partial^2_{\nu}u'\:dS.
\end{eqnarray*}
With the help of \eqref{upricc} and \eqref{lapbc1} we write
\begin{eqnarray*}
\partial^2_{\nu}u'=\Delta u'-(n-1)\partial_{\nu}u'\:H_{\partial\Omega}.
\end{eqnarray*}
Hence,
\begin{eqnarray*}
\int\limits_{\partial\Omega}(v\cdot\nu)\:\partial_{\nu}(\Delta u\:\partial_{\nu}u')\:dS
=
c_0\:\int\limits_{\partial\Omega}(v\cdot\nu)\:\Delta u'\:dS
-
c_0(n-1)\:\int\limits_{\partial\Omega}(v\cdot\nu)\:\partial_{\nu}u'\:H_{\partial\Omega}\:dS.
\end{eqnarray*}
Our computations yield a first simplification of \eqref{ddlamb}:
\begin{align*}
 \ddot{\Lambda}(0) = \int\limits_{\partial\Omega}\Delta u'\:\partial_\nu u'\:dS + \int\limits_{\partial\Omega}\Delta u\: \partial_\nu u''\:dS + c_0\int\limits_{\partial\Omega}(v\cdot\nu)\:\Delta u'\:dS.
\end{align*}
In the first integral on the right hand side we use \eqref{upricc} again. Thus, we get
\begin{eqnarray}\label{ddlamb1}
\ddot{\Lambda}(0)=c_0\:\int\limits_{\partial\Omega}\partial_{\nu}u''\:dS
\end{eqnarray}
In order to find a lower bound for $\ddot{\Lambda}(0)$, we analyze the integral in \eqref{ddlamb1}. Recall \eqref{tdot-eq5}. We differentiate this equation with respect to $t$ in $t=0$. Then \eqref{w-eq1} and \eqref{overrt} yield
\begin{eqnarray*}
&&\nu'\cdot\nabla u'
+
v\cdot D\nu\cdot\nabla u'
+
\partial_{\nu}u''
+
\nu\cdot D^2u'\cdot v
=\\
&&-
\Delta u'\:(v\cdot\nu)
-
c_0\:(v\cdot\nu')
-
c_0v\cdot D\nu\cdot v
-
c_0\:(w\cdot\nu).
\end{eqnarray*}
As before, $\nu'\cdot\nabla u'=0$ on $\partial\Omega$. Moreover, by \eqref{upricc} 
\begin{eqnarray*}
v\cdot D\nu\cdot\nabla u'=-
c_0v\cdot D\nu\cdot \nu\:(v\cdot\nu)=0,
\end{eqnarray*}
where the last equality follows from \eqref{normal}. Thus,
\begin{eqnarray}\label{lambdaddot2}
\ddot{\Lambda}(0)&=&
-
c_0\:\int\limits_{\partial\Omega}(v\cdot\nu)\:\Delta u'\:dS
-
c_0\:\int\limits_{\partial\Omega}\nu\cdot D^2u'\cdot v\:dS\\
\nonumber&&-
c_0^2\:\int\limits_{\partial\Omega}(v\cdot\nu')\:dS
-
c_0^2\:\int\limits_{\partial\Omega}v\cdot D\nu\cdot v\:dS
-
c_0^2\:\int\limits_{\partial\Omega}(w\cdot\nu)\:dS.
\end{eqnarray}
For the first integral we use \eqref{upricc} and we observe that Gau\ss\ theorem, partial integration and equation \eqref{uprimeq} for $u'$ gives
\begin{eqnarray}\label{int1}
-
c_0\:\int\limits_{\partial\Omega}(v\cdot\nu)\:\Delta u'\:dS
=
\int\limits_{\partial\Omega}\Delta u'\:\partial_{\nu}u'\:dS
=
\int\limits_{\Omega}\vert\Delta u'\vert^2\:dx-\Lambda(\Omega)\int\limits_{\Omega}\vert\nabla u'\vert^2\:dx.
\end{eqnarray}
The second intergal is slightly more involved. We set $v^{\tau}=v-(v\cdot\nu)\nu$. Since $\nabla u'=(\partial_{\nu}u')\nu$ and since \eqref{lapbc1} can be applied to $u'$, we get
\begin{eqnarray*}
-
c_0\:\int\limits_{\partial\Omega}v\cdot D^2u'\cdot\nu \:dS
&=&
-
c_0\:\int\limits_{\partial\Omega}v^{\tau}\cdot D^2u'\cdot \nu\:dS
-
c_0\:\int\limits_{\partial\Omega}(v\cdot\nu)\left(\Delta u'-(n-1)\partial_{\nu}u'\:H_{\partial\Omega}\right)\:dS\\
&=&
-
c_0\:\int\limits_{\partial\Omega}v^{\tau}\cdot D\left(\partial_{\nu}u'\:\nu\right)\cdot \nu\:dS
-
c_0\:\int\limits_{\partial\Omega}(v\cdot\nu)\:\Delta u'\:dS\\
&&-
c_0^2\:(n-1)\:\int\limits_{\partial\Omega}(v\cdot\nu)^2\:H_{\partial\Omega}\:dS.
\end{eqnarray*}
For the last equality we also used
\begin{eqnarray*}
v^{\tau}\cdot D\nu\cdot \nu=v^{\tau}\cdot D^{\tau}\nu\cdot \nu=0\quad\hbox{in}\:\partial\Omega.
\end{eqnarray*}
Next we note that with \eqref{upricc} we have
\begin{eqnarray*}
-
c_0\:\int\limits_{\partial\Omega}v^{\tau}\cdot D\left(\partial_{\nu}u'\:\nu\right)\cdot \nu\:dS
&=&
-
c_0\:\int\limits_{\partial\Omega}v^{\tau}\cdot D^{\tau}\left(\partial_{\nu}u'\:\nu\right)\cdot \nu\:dS
=
c_0^2\:\int\limits_{\partial\Omega}v^{\tau}\cdot D^{\tau}\left((v\cdot\nu)\:\nu\right)\cdot \nu\:dS\\
&=&
c_0^2\:\int\limits_{\partial\Omega}v^{\tau}\cdot \nabla^{\tau}(v\cdot\nu)\:dS,
\end{eqnarray*}
where the last equality uses \eqref{normal}.
\newline
For the third integral in \eqref{lambdaddot2} we apply formula \eqref{normal2}:
\begin{eqnarray*}
-
c_0^2\:\int\limits_{\partial\Omega}(v\cdot\nu')\:dS
=
c_0^2\:\int\limits_{\partial\Omega}\:v\cdot\nabla^{\tau}(v\cdot\nu)\:dS
=
c_0^2\:\int\limits_{\partial\Omega}\:v^{\tau}\cdot\nabla^{\tau}(v\cdot\nu)\:dS.
\end{eqnarray*}
These computations simplify \eqref{lambdaddot2} and we obtain
\begin{eqnarray}\label{lambdaddot3}
\ddot{\Lambda}(0)&=&
2\:\int\limits_{\partial\Omega}\partial_{\nu}u'\:\Delta u'\:dS
+
2c_0^2\:\int\limits_{\partial\Omega}v^{\tau}\cdot \nabla^{\tau}(v\cdot\nu)\:dS
-
c_0^2\:(n-1)\:\int\limits_{\partial\Omega}(v\cdot\nu)^2\:H_{\partial\Omega}\:dS\\
&&
\nonumber-
c_0^2\:\int\limits_{\partial\Omega}v\cdot D\nu\cdot v\:dS
-
c_0^2\:\int\limits_{\partial\Omega}(w\cdot\nu)\:dS.
\end{eqnarray}
Next we use the volume constraint \eqref{volume2}.
\begin{eqnarray*}
-
c_0^2\:\int\limits_{\partial\Omega}(w\cdot\nu)\:dS
&=&
c_0^2\:\int\limits_{\partial\Omega}(v\cdot\nu)\:\Div  v\:dS
-
c_0^2\:\int\limits_{\partial\Omega}  v\cdot Dv\cdot\nu\:dS\\
&=&
c_0^2\:\int\limits_{\partial\Omega}(v\cdot\nu)\:\Div_{\partial\Omega} v\:dS
-
c_0^2\:\int\limits_{\partial\Omega} v^{\tau}\cdot D^{\tau}v\cdot\nu\:dS.
\end{eqnarray*}
We integrate by parts in the first integral (see formula \eqref{tanpi} and \eqref{tgradient}).
\begin{eqnarray*}
-
c_0^2\:\int\limits_{\partial\Omega}(w\cdot\nu)\:dS
&=&
-c_0^2\:\int\limits_{\partial\Omega}v^{\tau}\cdot D^{\tau}(v\cdot\nu)\:dS
+
c_0^2(n-1)\:\int\limits_{\partial\Omega}(v\cdot\nu)^2\:H_{\partial\Omega}\:dS\\
&&-
c_0^2\:\int\limits_{\partial\Omega} v^{\tau}\cdot D^{\tau}v\cdot\nu\:dS.
\end{eqnarray*}
Thus, \eqref{lambdaddot3} becomes
\begin{eqnarray*}
\ddot{\Lambda}(0)&=&
2\:\int\limits_{\partial\Omega}\partial_{\nu}u'\:\Delta u'\:dS
+
c_0^2\:\int\limits_{\partial\Omega}v^{\tau}\cdot \nabla^{\tau}(v\cdot\nu)\:dS
-
c_0^2\:\int\limits_{\partial\Omega} v^{\tau}\cdot D^{\tau}v\cdot\nu\:dS\\
&&
\nonumber-
c_0^2\:\int\limits_{\partial\Omega}v\cdot D\nu\cdot v\:dS.
\end{eqnarray*}
An application of \eqref{normal} and  \eqref{normal2} yields
\begin{eqnarray*}
&&v^{\tau}\cdot \nabla^{\tau}(v\cdot\nu)
-
v^{\tau}\cdot D^{\tau}v\cdot\nu
-
v\cdot D\nu\cdot v
=
v^{\tau}\cdot D^{\tau}\nu\cdot v
-
v\cdot D\nu\cdot v\\
&&\qquad =
-(v\cdot\nu)\nu\cdot D\nu\cdot v=0.
\end{eqnarray*}
Thus, with \eqref{lambdaddot3} we proved the following lemma.
%%%%%%%%%%%%%%%%%%%%%%%%%%%%%%%%%%%%%%%%%
\begin{lemma}\label{lambdaddot4} Let $u'$ be the shape derivative of $u$ resulting from a volume preserving perturbation of $\Omega$. Then there holds
\[
\frac{\ddot{\Lambda}(0)}{2}= {\cal{E}}(u'),
\]
where
\[
{\cal{E}}(u') = \int\limits_{\Omega}\vert\Delta u'\vert^2\:dx-\Lambda(\Omega)\int\limits_{\Omega}\vert\nabla u'\vert^2\:dx.
\]
\end{lemma}
%%%%%%%%%%%%%%%%%%%%%%%%%%%%%%%%%%%%%%%%%
%%%%%%%%%%%%%%%%%%%%%%%%%%%%%%%%%%%%%%%%%
\section{Minimization of the second domain variation }
%%%%%%%%%%%%%%%%%%%%%%%%%%%%%%%%%%%%%%%%%
%%%%%%%%%%%%%%%%%%%%%%%%%%%%%%%%%%%%%%%%%
In this chapter we consider the quadratic functional 
\begin{eqnarray}\label{E1}
  {\cal{E}}(\varphi) := \int\limits_{\Omega}\vert \Delta \varphi\vert^2\:dx - \Lambda(\Omega)\int\limits_{\Omega}\vert\nabla \varphi\vert^2\:dx
\end{eqnarray}
for $\varphi \in H^{1,2}_0\cap H^{2,2}(\Omega)$. It will be convenient to work with an alternative representation of $ {\cal{E}}$. For $\varphi\in H^{1,2}_0\cap H^{2,2}(\Omega)$  there holds
\begin{eqnarray*}
{\cal{E}}(\varphi) = \int\limits_\Omega \vert D^2\varphi\vert^2 - \Lambda(\Omega)\vert\nabla\varphi\vert^2\:dx + \int\limits_{\partial\Omega}\Delta\varphi\partial_\nu\varphi - \varphi\cdot D^2\varphi\cdot\nu\:dS.
\end{eqnarray*}
We apply \eqref{lapbc1} and \eqref{normal}.
\begin{eqnarray*}
  \Delta\varphi\partial_\nu\varphi - \varphi\cdot D^2\varphi\cdot\nu &=& \partial^2_\nu\varphi\,\partial_\nu\varphi + (n-1)(\partial_\nu\varphi)^2\,H_{\partial\Omega} - \varphi\cdot D^2\varphi\cdot\nu \\
 &=& \nu\cdot D^2\varphi\cdot\nu \,(\nu\cdot\nabla\varphi) + (n-1)(\partial_\nu\varphi)^2\,H_{\partial\Omega} - \varphi\cdot D^2\varphi\cdot\nu \\
&=& (n-1)(\partial_\nu\varphi)^2\,H_{\partial\Omega}.
\end{eqnarray*}
Consequently, we get
\begin{eqnarray}\label{E2}
  {\cal{E}}(\varphi) = \int\limits_{\Omega}\vert D^2\varphi\vert^2\:dx - \Lambda(\Omega)\int\limits_{\Omega}\vert\nabla \varphi\vert^2\:dx + (n-1)\int\limits_{\partial\Omega}(\partial_\nu \varphi)^2H_{\partial\Omega}\:dS.
\end{eqnarray}
%%%%%%%%%%%%%%%%%%%%%%%%%%%%%%%%%%%%%%%%%%%%%%
\begin{remark}\label{lsc}
 The functional $ {\cal{E}}$ is lower semicontinuous with respect to weak convergence in $H^{1,2}_0\cap H^{2,2}(\Omega)$. 
\end{remark}
%%%%%%%%%%%%%%%%%%%%%%%%%%%%%%%%%%%%%%%
Since $\Omega$ is optimal, we know from Lemma \ref{lambdaddot4} that 
\[
  {\cal{E}}(\varphi)\geq 0
\]
for all $\varphi$ which are shape derivatives of $u$. Recall that $\varphi$ is a shape derivative, if it solves \eqref{shape1} - \eqref{shape4} for some vector field $v$ in the class described in Chapter 1 (Lemma \ref{volpres1}).
\newline
\newline 
The following remark shows a property of shape derivatives we have not yet mentioned.
%%%%%%%%%%%%%%%%%%%%%%%%%%%%%%%%%%%%%%%%%%%%%%
\begin{remark}\label{d_nuprime}
Let $\varphi$ be a shape derivative and assume that $\partial_\nu \varphi\equiv0$ in $\partial\Omega$. Then $\varphi \in H^{2,2}_0(\Omega)$ and, since $\varphi$ satisfies equation \eqref{uprimeq}, $\varphi$ is a buckling eigenfunction in $\Omega$. Thus by uniqueness of $u$ we get $\varphi = \alpha u$ for any $\alpha \in\mathbb{R}$. Then formula \eqref{Rel} yields
\begin{eqnarray*}
\Lambda(\Omega) =  \int\limits_{\partial\Omega}\vert\Delta \varphi\vert^2\, x\cdot\nu\:dS = \alpha^2c_o^2\int\limits_{\partial\Omega}x\cdot\nu\:dS =  \alpha^2\int\limits_{\partial\Omega}\vert\Delta u\vert^2\, x\cdot\nu\:dS =\alpha^2 \Lambda(\Omega).
\end{eqnarray*}
Thus, $\alpha^2=1$ and there holds
\[
\left\vert \int\limits_\Omega \nabla u\cdot\nabla \varphi \:dx \right\vert =1.
\]
This is contradictory to \eqref{graduprime} and thus $\partial_{\nu} \varphi$ cannot vanish identically on $\partial\Omega$.
\end{remark}
%%%%%%%%%%%%%%%%%%%%%%%%%%%%%%%%%%%%%%%%%
This motivates the following definition.
\begin{eqnarray*}
{\cal{Z}}:= \left\{ \varphi \in H^{1,2}_0\cap H^{2,2}(\Omega) : \int\limits_{\partial\Omega}\partial_\nu\varphi\:dS=0,\,  \int\limits_{\partial\Omega}(\partial_\nu\varphi)^2\:dS>0,\,  \int\limits_{\Omega}\nabla u\cdot \nabla\varphi\:dx=0 \right\}.
\end{eqnarray*}
Note that ${\cal{Z}}$ contains elements which are not shape derivatives. Nevertheless we will show that 
\[
   {\cal{E}}\big{|}_{{\cal{Z}}}\geq 0.
\] 
%%%%%%%%%%%%%%%%%%%%%%%%%%%%%%%%%%%%%%%%%
The next lemma ensures that ${\cal{Z}}$ is not empty and that at least for a specific shape derivative $\cal{E}$ is equal to zero. 
\begin{lemma}\label{E(dku)}
 For each $1\leq k\leq n$ the directional derivative $\partial_ku$ satisfies $\partial_ku\in \cal{{\cal{Z}}}$. Furthermore, ${\cal{E}}(\partial_ku)=0$.
\end{lemma}
{\bf{Proof}} Let $1 \leq k\leq n$. Due to \eqref{eq1} and \eqref{eq2} $\partial_ku$ satisfies
\begin{eqnarray}
 \Delta^2\partial_ku+\Lambda(\Omega)\Delta \partial_ku&=&0\quad\hbox{in}\:\Omega \label{d_ku_pde}\\
\partial_k u&=&0\quad\hbox{in}\:\partial\Omega. \notag
\end{eqnarray}
According to \eqref{lapbc2} there holds $\partial_\nu\partial_ku = c_0\nu_k$ on $\partial\Omega$. Hence, 
\[
  \int\limits_{\partial\Omega}\partial_\nu\partial_ku\:dS =  c_0 \int\limits_{\partial\Omega}\nu_k\:dS =0.
\]
In addition, we find that
\[
  \int\limits_\Omega \nabla u\cdot\nabla\partial_ku\;dx = \frac{1}{2}\int\limits_{\partial\Omega}\vert\nabla u\vert^2\nu_k\:dS =0.
\]
Following the idea of Remark \ref{d_nuprime}, we obtain that $\partial_\nu\partial_ku$ does not vanish identically on $\partial\Omega$. Thus, $\partial_ku \in \cal{{\cal{Z}}}$. 
Moreover, \eqref{w-eq1} and \eqref{d_ku_pde} imply
\[
  {\cal{E}}(\partial_ku) = \int\limits_\Omega (\Delta^2\partial_ku + \Lambda( \Omega)\Delta\partial_ku)\partial_ku\:dx + \int\limits_{\partial\Omega}\partial_k\Delta u\, \partial_\nu\partial_ku\:dS =0.
\]
This proves the lemma.
\hfill$\square$
\newline
\newline
Note that each directional derivative of $u$ is a shape derivative resulting from translations of $\Omega$.
%%%%%%%%%%%%%%%%%%%%%%%%%%%%%%%%%%%%%%%%%
\begin{theorem} \label{inffinite}
 The infimum of the functional $\cal{E}$ in ${\cal{Z}}$ is finite. 
\end{theorem}
%%%%%%%%%%%%%%%%%%%%%%%%%%%%%%%%%%%%%%%%%
{\bf{Proof}} We argue by contradiction. Let us assume that $\inf_{\cal{{\cal{Z}}}} {\cal{E}} = -\infty$ and consider a sequence $(\hat{w}_k)_k \subset \cal{{\cal{Z}}}$ such that 
\[
  \lim_{k\to\infty} {\cal{E}}(\hat{w}_k) = -\infty.
\]
For this sequence there either holds
\[
  \int\limits_{\partial\Omega}(\partial_\nu\hat{w}_k)^2\:dS \stackrel{k\to\infty}{\longrightarrow}0  \qquad
\hbox{or} \qquad \int\limits_{\partial\Omega}(\partial_\nu\hat{w}_k)^2\:dS \stackrel{k\to\infty}{\not\longrightarrow}0.
\]
If the second case holds true, we normalize the sequence $(\hat{w}_k)_k$ such that $\vert\vert\partial_\nu\hat{w}_k\vert\vert_{L^2(\partial\Omega)}=1$. Hence, in either case, for each $\hat{w}_k$ there holds $\|\partial_\nu\hat{w}_k\|_{L^2(\partial\Omega)}\leq1$. Thus, \eqref{H_dOmega} gives
\begin{eqnarray*}
\left\vert\, \int\limits_{\partial\Omega} H_{\partial\Omega}\, (\partial_\nu\hat{w}_k)^2\:dS\right\vert \leq \max_{\partial\Omega}\vert H_{\partial\Omega}\vert < \infty.
\end{eqnarray*}
We use \eqref{E2} and obtain
\begin{eqnarray}\label{estE1}
{\cal{E}}(\hat{w}_k) \geq -\Lambda(0) \int\limits_\Omega \vert\nabla\hat{w}_k\vert\:dx 
- (n-1) \max_{\partial\Omega}\vert H_{\partial\Omega}\vert.
\end{eqnarray}
The assumption $\lim_{k\to \infty}{\cal{E}}(\hat{w_k}) = -\infty$ implies 
\begin{eqnarray*}
   \int\limits_\Omega \vert\nabla\hat{w}_k\vert^2\:dx \stackrel{k\to\infty}{\longrightarrow} \infty.
\end{eqnarray*}
We define 
\begin{eqnarray*}
   w_k := \frac{1}{\vert\vert\nabla\hat{w}_k\vert\vert_{L^2(\Omega)}}\hat{w}_k. 
\end{eqnarray*}
Then there holds
\begin{eqnarray}\label{w_k}
  \vert\vert\nabla w_k\vert\vert_{L^2(\Omega)} = 1 \quad \hbox{and}\quad \int\limits_{\partial\Omega}(\partial_\nu w_k)^2\:dS\;\stackrel{k\to\infty}{\longrightarrow}\; 0.
\end{eqnarray}
Moreover, for each $k\in\mathbb{N}$ estimate \eqref{estE1} implies
\[
   {\cal{E}}(w_k) \geq -\Lambda(0) - C
\]
and the infimum of $\cal{E}$ in $M := \{w_k : k\in\mathbb{N}\}$ is finite. Therefore, we can choose a subsequence of $(w_k)_k$, denote by $(w_k)_k$ as well, such that 
\[
   \lim_{k\to\infty} {\cal{E}}(w_k) = \inf_M \cal{E}.
\]
Now Poincar\'e's inequality and the previous estimates imply
\begin{eqnarray*}
  \vert\vert w_k\vert\vert^2_{H^{2,2}(\Omega)} &=& \int\limits_\Omega \vert D^2w_k\vert^2 + \vert\nabla w_k\vert^2 + w_k^2\:dx \\
  &\leq& {\cal{E}}(w_k) + C \int\limits_\Omega\vert\nabla w_k\vert^2\:dx + (n-1)\int\limits_{\partial\Omega}\vert H_{\partial\Omega}\vert(\partial_\nu w_k)^2\:dS \\
 &\leq& C.
\end{eqnarray*}
Thus, the sequence $(w_k)_k$ is uniformly bounded in $H^{2,2}(\Omega)$ and there exists a $w \in H^{2,2}(\Omega)$ such that $(w_k)_k$ weakly converges to $w$. In view of \eqref{w_k}, the limit function $w$ satisfies $ \vert\vert\nabla w\vert\vert_{L^2(\Omega)} = 1 $ and $\partial_\nu w =0$ on $\partial\Omega$. Since $w_k =0$ in $\partial\Omega$ for each $k \in \mathbb{N}$, we conclude that $w \in H^{2,2}_0(\Omega)$. 
\newline
Now let us recall that ${\cal{E}}(\hat{w}_k)$ converges to $-\infty$. Thus there exists a $k_0 \in \mathbb{N}$ such that 
\begin{eqnarray*}
{\cal{E}}(w_k) = \frac{1}{\vert\vert\nabla\hat{w}_k\vert\vert_{L^2(\Omega)}} \,{\cal{E}}(\hat{w}_k) < 0
\end{eqnarray*}
for all $k\geq k_0$. Since the functional $\cal{E}$ is lower semicontinous with respect to weak convergence in $H^{2,2}(\Omega)$, we find that ${\cal{E}}(w)<0$. According to the definiton of $\cal{E}$ in \eqref{E1}, this immediately leads to
\begin{eqnarray*}
  \frac{\int\limits_{\Omega}\vert\Delta w\vert^2\:dx}{\int\limits_\Omega\vert\nabla w\vert^2\:dx} < \Lambda(\Omega).
\end{eqnarray*}
Since $w \in H^{2,2}_0(\Omega)$ this is contradictory to the minimum property of $\Lambda(\Omega)$. 
\hfill $\square$
\newline
\newline
%%%%%%%%%%%%%%%%%%%%%%%%%%%%%%%%%%%%%%%%%%%%%%%
As mentioned in the previous proof, a minimizing sequence $(\varphi_k)_k\subset {\cal{Z}}$ satisfies one of the following two conditions
\begin{itemize}
 \item [i)] $ \int\limits_{\partial\Omega}(\partial_\nu \varphi_k)^2\:dS \stackrel{k\to\infty}{\longrightarrow} 0$
\item[ii)]$\int\limits_{\partial\Omega}(\partial_\nu\varphi_k)^2\:dS \stackrel{k\to\infty}{\not\longrightarrow}0$.
\end{itemize}
In the sequel, we show that the case i) implies that the minimizing sequence $(\varphi_k)_k$ converges to zero. For this purpose, let $(\varphi_k)_k\subset {\cal{Z}}$ be a minimizing sequence which satisfies condition i). 
From \eqref{E2} we get
\begin{eqnarray*}
\vert\vert\varphi_k\vert\vert^2_{H^{2,2}(\Omega)} \leq C,
\end{eqnarray*}
and thus there exists a $\varphi\in H^{2,2}(\Omega)$ such that $\varphi_k$ weakly converges to $\varphi$ in $H^{2,2}(\Omega)$ and ${\cal{E}}(\varphi)= \inf_{\cal{Z}} {\cal{E}}$.  Furthermore, condition i) implies $\varphi \in H^{2,2}_0(\Omega)$. From Lemma \ref{E(dku)} we obtain 
\begin{eqnarray*}
\inf_{\cal{Z}}{\cal{E}}={\cal{E}}(\varphi) \leq {\cal{E}}(\partial_lu) =0 \quad \hbox{for any } 1\leq l \leq n.
\end{eqnarray*}
Hence,
\begin{eqnarray*}
\frac{\int\limits_\Omega\vert\Delta\varphi\vert^2\:dx}{\int\limits_\Omega \vert\nabla\varphi\vert^2\:dx} \leq \Lambda(\Omega).
\end{eqnarray*}
Thus $\varphi$ is necessarily an eigenfunction corresponding to $\Lambda(\Omega)$. Since the eigenvalue is simple we have $\varphi = \alpha\,u$ for $\alpha \in \mathbb{R}$. Now let us recall that $\varphi_k \in {\cal{Z}}$ and, therefore, 
\begin{eqnarray*}
  \alpha =  \alpha\,\int\limits_\Omega \vert\nabla u\vert^2\:dx =   \int\limits_\Omega \nabla u\cdot\nabla\varphi\:dx =  \lim_{k\to\infty} \int\limits_\Omega \nabla u\cdot\nabla\varphi_k\:dx =0. 
\end{eqnarray*} 
Consequently, $\alpha =0$ and $\varphi \equiv 0$ in $\Omega$. Hence $\varphi \notin {\cal{Z}}$. Since we are interested to find minimizers of ${\cal{E}}$ in ${\cal{Z}}$, we restrict ourselves to minimzing sequences which satisfy the condition ii). Thus we consider the functional
\begin{eqnarray*}
   \tilde{{\cal{E}}}(\varphi) := \frac{{\cal{E}}(\varphi)}{\int\limits_{\partial\Omega}(\partial_\nu\varphi)^2\:dS},
\end{eqnarray*}
where $\varphi \in {\cal{Z}}$ and we set $\tilde{{\cal{E}}}=\infty$ if $\int\limits_{\partial\Omega}(\partial_\nu\varphi_k)^2\:dS=0$.
%%%%%%%%%%%%%%%%%%%%%%%%%%%%%%%%%%%%%%%%%%%%%%
\begin{remark}\label{gradmise}
 Suppose $(\varphi_k)_k \subset {\cal{Z}}$ is a minimzing sequence for $\tilde{{\cal{E}}}$ in ${\cal{Z}}$. Then there exists a constant $C>0$ such that $\vert\vert\nabla\varphi_k\vert\vert_{L^2(\Omega)}\leq C$ for every $k\in\mathbb{N}$. This follows by contradiction.
Otherwise we may assume that $\vert\vert\nabla\varphi_k\vert\vert_{L^2(\Omega)}$ tends to infinity as $k\to\infty$, we define $\varphi^\ast_k := \vert\vert\nabla\varphi_k\vert\vert_{L^2(\Omega)}^{-1}\varphi_k$. Then $(\varphi^\ast_k)_k$ is uniformly bounded in $H^{2,2}(\Omega)$ and 
\begin{eqnarray*}
 \int\limits_{\partial\Omega}(\partial_\nu \varphi^\ast_k)^2\:dS \stackrel{k\to\infty}{\longrightarrow} 0.
\end{eqnarray*}
Thus, $(\varphi^\ast_k)_k$ converges weakly to a function $\varphi \in H^{2,2}_0(\Omega)$ and for every $1 \leq l\leq n$ there holds 
\begin{eqnarray*}
  \inf_{\cal{Z}} \tilde{{\cal{E}}} = \tilde{{\cal{E}}}(\varphi) \leq \tilde{{\cal{E}}}(\partial_lu) =0.
\end{eqnarray*}
As the previous considerations have shown, this implies $\varphi \equiv 0$ in $\Omega$. Thus, our assumption cannot be true. 
\end{remark}
%%%%%%%%%%%%%%%%%%%%%%%%%%%%%%%%%%%%%%%%
We now consider a minimizing sequence $(\varphi_k)_k \subset {\cal{Z}}$ which satisfies 
\begin{eqnarray}\label{normmise}
\int\limits_{\partial\Omega}(\partial_\nu\varphi_k)^2\:dS=1
\end{eqnarray}
for all $k\in\N$. As before we obtain the inequality
\begin{eqnarray*}
\vert\vert\varphi_k\vert\vert^2_{H^{2,2}(\Omega)} \leq {\cal{E}}(\varphi_k) + C\int\limits_\Omega\vert \nabla\varphi_k\vert^2\:dx .
\end{eqnarray*}
Thus, $(\varphi_k)_k$ is uniformly bounded in $H^{2,2}(\Omega)$  and $\varphi_k$ converges weakly to a $\varphi^\ast \in H^{2,2}(\Omega)$. We find that $\varphi^\ast \in {\cal{Z}}$ and 
${\cal{E}}(\varphi^\ast) = \inf_{\cal{Z}}E$. In addition, there holds
\[
   \int\limits_{\partial\Omega}(\partial_\nu\varphi^\ast)^2\:dS = 1. 
\]
Hence,  $\varphi^\ast$ minimizes $\tilde{{\cal{E}}}$ in ${\cal{Z}}$. Suppose $\theta \in {\cal{Z}}$, then the minimality of $\varphi^\ast$ implies
\[
 \frac{d}{dt}\frac{{\cal{E}}(\varphi^\ast +t \theta)}{\int\limits_{\partial\Omega}(\partial_\nu (\varphi^\ast + t \theta))^2\:dS)}\bigg\vert_{t=0} =0
\]
and we obtain
\[
 \int\limits_\Omega [\Delta^2\varphi^\ast + \Lambda(\Omega)\Delta\varphi]\,\theta\:dx - \int\limits_{\partial\Omega}[\Delta\varphi^\ast + \rho\,\partial_\nu\varphi^\ast]\,\partial_\nu\theta\:dS =0. 
\]
Since $\theta\in {\cal{Z}}$ was chosen arbitrary, $\varphi^\ast$ satisfies the Euler-Lagrange equalities
\begin{eqnarray*}
 \Delta^2\varphi^\ast + \Lambda(\Omega)\Delta\varphi^\ast &=& 0 \quad \hbox{in} \quad \Omega \\
 \Delta\varphi^\ast + \rho\partial_\nu\varphi^\ast &=& const. \quad \hbox{in} \quad \partial\Omega,
\end{eqnarray*}
where $\rho := \min_{\cal{Z}}\tilde{{\cal{E}}}$. The following theorem collects the previous results.
%%%%%%%%%%%%%%%%%%%%%%%%%%%%%%%%%%%%%%%%
\begin{theorem}
 There exists a function $\varphi^\ast \in {\cal{Z}}$ such that $\tilde{{\cal{E}}}(\varphi^\ast)=\min_{\cal{Z}}\tilde{{\cal{E}}}$. Furthermore, any minimizer $\varphi^\ast \in {\cal{Z}}$ satisfies
\begin{eqnarray}
\label{EL1} \Delta^2\varphi^\ast + \Lambda(\Omega)\Delta\varphi^\ast &=& 0 \quad \hbox{in} \quad \Omega \\
\label{EL2} \Delta\varphi^\ast + \rho\,\partial_\nu\varphi^\ast &=& const. \quad \hbox{in} \quad \partial\Omega\\
\notag\varphi^\ast &=&0 \quad \hbox{in} \quad \partial\Omega,
\end{eqnarray}
where $\rho := \min_{\cal{Z}}\tilde{{\cal{E}}}$.
\end{theorem}
%%%%%%%%%%%%%%%%%%%%%%%%%%%%%%%%%%%%%%%%%%%%%
The next theorem shows that in fact $\rho =0$.
%%%%%%%%%%%%%%%%%%%%%%%%%%%%%%%%%%%%%%%%%%%%%
\begin{theorem}\label{minE=0}
Suppose $\varphi^{\ast}\in {\cal{Z}}$ is a minimizer of $\tilde{{\cal{E}}}$. Then there holds $\tilde{{\cal{E}}}(\varphi^\ast) =0$. In particular, ${\cal{E}}\geq 0$ in $\cal{Z}$.
\end{theorem}
%%%%%%%%%%%%%%%%%%%%%%%%%%%%%%%%%%%%%%%%%%%%%
{\bf{Proof}}  Let $\varphi^\ast \in {\cal{Z}}$ be a minimizer of $\tilde{{\cal{E}}}$. Since $\varphi^\ast$ satisfies equation \eqref{EL1} and $\partial\Omega$ is smooth, $\varphi^\ast$ is a smooth function on $\overline{\Omega}$. Hence, we may define a volume preserving perturbation $\Phi_t$ of $\Omega$ such that 
\begin{eqnarray*}
  \partial_\nu u'(x) = \partial_\nu \varphi^\ast(x) \quad \hbox{for} \quad x \in   \partial\Omega.
\end{eqnarray*}
Note that this can be achieved by setting $v = c_0^{-1}\nabla \varphi^\ast$ in $\partial\Omega$. In this way, each minimizer $\varphi^{\ast}$ implies the existence of vector fields $v$ and $w$ in the sense of Section \ref{DomVar}. We define $\psi := u' - \varphi^\ast$, then $\psi \in H^{2,2}_0(\Omega)$ and
\begin{eqnarray*}
\Delta^2\psi + \Lambda(\Omega)\Delta\psi =0 \quad \hbox{in} \quad \Omega.
\end{eqnarray*}
The uniqueness of $u$ implies $\psi = \alpha \,u$ for an $\alpha \in \mathbb{R}$. Since $\varphi^\ast\in {\cal{Z}}$, equation \eqref{graduprime} yields
\begin{align*}
  0 = \int\limits_\Omega \nabla u\cdot\nabla u'\:dx -  \int\limits_\Omega \nabla u\cdot\nabla \varphi^\ast\:dx =  \int\limits_\Omega \nabla u\cdot\nabla \psi\:dx = \alpha.
\end{align*} 
Consequently, $u' \equiv \varphi^\ast$. Thus $\varphi^\ast$ is a shape derivative. Since $\Omega$ is optimal $\tilde{{\cal{E}}}(\varphi^\ast)\geq 0$. Finally we apply Lemma \ref{E(dku)}. This gives
\begin{eqnarray*}
   0\leq \tilde{{\cal{E}}}(\varphi^\ast) = \min_{\cal{Z}} \tilde{{\cal{E}}} \leq \tilde{{\cal{E}}}(\partial_ku) =0.
\end{eqnarray*}
\hfill $\square$
%%%%%%%%%%%%%%%%%%%%%%%%%%%%%%%%%%%%%%%%%%%
%%%%%%%%%%%%%%%%%%%%%%%%%%%%%%%%%%%%%%%%%%%
\section{The optimal domain is a ball}
%%%%%%%%%%%%%%%%%%%%%%%%%%%%%%%%%%%%%%%%%%%
%%%%%%%%%%%%%%%%%%%%%%%%%%%%%%%%%%%%%%%%%%%
We will use an inequaltiy due to L.E. Payne to show that the optimal domain $\Omega$ is a ball. Payne's inequality (see \cite{P55}) states that for each domain $G$ there holds 
\begin{equation*}
\lambda_2(G) \leq \Lambda(G)
\end{equation*}
and equality only holds if and only if $G$ is a ball. Thereby $\lambda_2$ denotes the second Dirichlet eigenfunction of the Laplacian.
In the sequel, we construct a suitable function $\psi \in {\cal{Z}}$ such that the condition ${\cal{E}}(\psi)\geq 0$ (due to Theorem \ref{minE=0}) will imply that the optimal domain $\Omega$ is a ball.
For this purpose, we denote by $u_1$ and $u_2$ the first and the second Dirichlet eigenfunction for the Laplacian in $\Omega$. Thus, for $1\leq k\leq 2$ there holds
\begin{eqnarray*}
\Delta u_k + \lambda_k(\Omega) u_k &=&0 \quad \mbox{in }\Omega \\
u_k &=&0 \quad \mbox{in }\partial\Omega,
\end{eqnarray*}
where $\lambda_k(\Omega)$ is the $k$-th Dirichlet eigenvalue for the Laplacian in $\Omega$. Note that $0< \lambda_1(\Omega)<\lambda_2(\Omega)$. For the sake of brevity, we will write $ \lambda_k$ instead of $\lambda_k(\Omega)$ and $\Lambda$ instead  of $\Lambda(\Omega)$. In addition, we assume $\vert\vert u_k\vert\vert_{L^2(\Omega)}=1$ and 
\begin{eqnarray*}
  \int\limits_\Omega u_1u_2\:dx =0.
\end{eqnarray*}
Without loss of generality, we may assume that 
\begin{eqnarray*}
   \int\limits_\Omega u_1\:dx >0 \quad \mbox{and} \quad \int\limits_\Omega u_2\:dx \leq 0.
\end{eqnarray*}
Consequently, there exists a $t \in (0,1]$ such that 
\begin{eqnarray}\label{meanval}
\int\limits_\Omega (1-t)\:\lambda_1\:u_1 + t\:\lambda_2\: u_2 \:dx =0.
\end{eqnarray}
This fixes $t$. Next we define 
\begin{eqnarray*}
\psi(x):= (1-t)\:u_1(x) + t\:u_2(x) + c\:u(x)\qquad\hbox{for}\:x\in\overline{\Omega},
\end{eqnarray*} 
where $u$ is the first buckling eigenfunction in $\Omega$. The constant $c$ is given by
\begin{eqnarray*}
c := -\frac{1}{\Lambda}\int\limits_\Omega (1-t)\lambda_1\nabla u.\nabla u_1 + t\lambda_2\nabla u.\nabla u_2\: dx.
\end{eqnarray*} 
In a first step we show that $\psi  \in {\cal{Z}}$. Note that $\psi \in H^{1,2}_0\cap H^{2,2}(\Omega)$. Moreover the definition of $\psi$, the fact that $\partial_{\nu}u=0$ on $\partial\Omega$, the equations for $u_1$ and $u_2$,   and \eqref{meanval} imply
\begin{eqnarray*}
\int\limits_{\partial\Omega} \partial_\nu \psi  \:dS 
=
\int\limits_\Omega (1-t)\:\Delta u_1 + t\:\Delta u_2\:dx
= 
-\int\limits_\Omega (1-t)\lambda_1u_1 + t\lambda_2u_2\:dx =0.
\end{eqnarray*} 
By the unique continuation principle $\partial_\nu \psi $ does not vanish identically in $\partial\Omega$. Thus, to show that $\psi  \in {\cal{Z}}$, it remains to prove that 
\begin{equation}\label{cond_3}
  \int\limits_\Omega \nabla u.\nabla \psi \:dx =0.
\end{equation}
We recall that $\Delta u = c_0$ in $\partial\Omega$. Hence 
\begin{align*}
 0 &= \int\limits_\Omega (\Delta^2u + \Lambda \Delta u)\psi  \:dx 
 = \int\limits_\Omega\Delta u \Delta \psi \:dx - \Lambda\int\limits_{\Omega}\nabla u.\nabla \psi \:dx \\
&= -\int\limits_\Omega [(1-t)\lambda_1u_1 + t\lambda_2u_2]\Delta u\:dx + c\int\limits_\Omega\vert\Delta u\vert^2\:dx- \Lambda\int\limits_\Omega\nabla u.\nabla \psi \:dx.
\end{align*}
Since $\vert\vert \nabla u\vert\vert_{L^2(\Omega)}=1$, the second integral is equal to $\Lambda$. Thus, the definition of $c$ implies \eqref{cond_3}. Note that $\psi$ is not a shape derivative since it fails to satisfy \eqref{uprimeq} - unless $t=1$ and $\Omega$ equals a ball.  However, $\psi  \in {\cal{Z}}$ and, according to Theorem \ref{minE=0}, there holds $\tilde{{\cal{E}}}(\psi)\geq 0$. Consequently, ${\cal{E}}(\psi)\geq 0$. Thus
\begin{align*}
  {\cal{E}}(\psi) &= \int\limits_\Omega \vert\Delta \psi\vert^2-\Lambda \vert\nabla \psi \vert^2\:dx  \\
& = (1-t)^2\lambda_1 (\lambda_1-\Lambda) + t^2\lambda_2(\lambda _2- \Lambda) + 2\,c\,c_0 \int\limits_\Omega (1-t)\lambda_1u_1 + t\lambda_2u_2\:dx \\
&\stackrel{\eqref{meanval}}{=} (1-t)^2\lambda_1 (\lambda_1-\Lambda) + t^2\lambda_2(\lambda _2- \Lambda) \geq 0.
\end{align*}
Since $\lambda_1-\Lambda < 0$ and $\lambda_2-\Lambda\leq 0$, both summands in ${\cal{E}}(\psi)$ have to vanish. Consequently $t=1$ and $\lambda_2(\Omega) = \Lambda(\Omega)$. Payne's inequality implies that $\Omega$ is a ball. This proves the main theorem of the paper.
%%%%%%%%%%%%%%%%%%%%%%%%%%%%%%%%%%%%%%%%%%
\begin{theorem}
Let $\Omega$ be a bounded, smooth and simply connected domain in $\R$, which minimizes the first buckling eigenvalue among all bounded, smooth and simply connected domains in $\R$ with given measure. Then $\Omega$ is a ball.
\end{theorem} 
%%%%%%%%%%%%%%%%%%%%%%%%%%%%%%%%%%%%%%%%%
%%%%%%%%%%%%%%%%%%%%%%%%%%%%%%%%%%%%%%%%%
%{\bf Acknowledgement} 
%%%%%%%%%%%%%%%%%%%%%%%%%%%%%%%%%%%%%%%%%%%%%
%%%%%%%%%%%%%%%%%%%%%%%%%%%%%%%%%%%%%%%%%%%%%

\end{document}